\documentclass[12pt, reqno]{amsart}
\usepackage{amssymb,amsmath}
\usepackage{graphicx}
\usepackage{hyperref}
\hypersetup{colorlinks=true,linkcolor=black,anchorcolor=black,citecolor=black}

\usepackage{appendix}
\usepackage{xcolor}
\usepackage{longtable}
\usepackage{float}
\usepackage{tikz}
\usepackage[all]{xy}

\DeclareFontFamily{U}{mathb}{\hyphenchar\font45}
\DeclareFontShape{U}{mathb}{m}{n}{
      <5> <6> <7> <8> <9> <10> gen * mathb
      <10.95> mathb10 <12> <14.4> <17.28> <20.74> <24.88> mathb12
      }{}
\DeclareSymbolFont{mathb}{U}{mathb}{m}{n}
\DeclareMathSymbol{\righttoleftarrow}{3}{mathb}{"FD}

\theoremstyle{plain}
\newtheorem{prop}{Proposition}[section]

\newtheorem{theo}[prop]{Theorem}
\newtheorem{coro}[prop]{Corollary}
\newtheorem{lemm}[prop]{Lemma}
\theoremstyle{remark}

\theoremstyle{definition}
\newtheorem{rema}[prop]{Remark}

\newtheorem{exam}[prop]{Example}
\numberwithin{equation}{section}

\newcommand{\bP}{\mathbf P}

\def\cI{{\mathcal I}}

\def\cN{{\mathcal N}}

\def\cP{{\mathcal P}}
\def\cQ{{\mathcal Q}}

\def\cX{{\mathcal X}}

\def\cZ{{\mathcal Z}}

\def\wt{{\mathrm {wt}}}

\newcommand{\Pf}[1]{\langle\!\langle #1\rangle\!\rangle}

\def\bA{{\mathbb A}}

\def\bP{{\mathbb P}}

\def\bZ{{\mathbb Z}}

\def\bC{{\mathbb C}}

\newcommand{\nr}{\mathrm{nr}}
\def\rH{{\mathrm H}}

\def\CH{\mathrm{CH}}

\def\lim{\mathrm{lim}}

\def\Spec{\mathrm{Spec}}

\def\Sing{\mathrm{Sing}}

\title[Quartic hypersurfaces]{Stable irrationality of quartic hypersurfaces}

\author[B\"ohning]{Christian B\"ohning}\thanks{For the purpose of open access, the authors have applied a Creative Commons Attribution (CC BY) licence to any Author Accepted Manuscript version arising from this submission.}
\address{Christian B\"ohning, Mathematics Institute, University of Warwick\\
Coventry CV4 7AL, England}
\email{C.Boehning@warwick.ac.uk}

\author[von Bothmer]{Hans-Christian Graf von Bothmer}
\address{Hans-Christian Graf von Bothmer, Fachbereich Mathematik der Universit\"at Hamburg\\
Bundesstra\ss e 55\\
20146 Hamburg, Germany}
\email{hans.christian.v.bothmer@uni-hamburg.de}

\author[Pirutka]{Alena Pirutka}
\address{Alena Pirutka, Courant Institute of Mathematical Sciences, New York University\\
251 Mercer Street\\
NY 10012, USA}
\email{ap166@nyu.edu}

\author[Zhang]{Zhijia Zhang}

\address{Zhijia Zhang, 
Department of Mathematics, Columbia University, New York, NY 10027, USA
}

\email{zhang.zhijia@columbia.edu}

\date{\today}

\makeatletter
\@namedef{subjclassname@2020}{\textup{2020} Mathematics Subject Classification}
\makeatother

\begin{document}

\begin{abstract}
We prove that a very general quartic hypersurface of any positive dimension over an uncountable algebraically closed field of characteristic different from 2 does not admit an integral decomposition of the diagonal, and in particular it is not stably rational.
\end{abstract}
	
	\maketitle

    \section{Introduction}
	\label{sect:intro}

A variety $X$ over a field $k$ is called {\em rational} if $X$ is birational to $\bP^n$ over $k$, and it is called {\em stably rational} if $X\times\bP^r$ is rational for some integer $r\geq0$. A classical problem in algebraic geometry is to determine whether a given variety is (stably) rational, usually referred to as the {\em rationality problem}. One of the most fundamental cases is the rationality problem for hypersurfaces, that is, varieties in $\bP^n$ defined by a single homogeneous polynomial.

The rationality problem for hypersurfaces has a long history and has drawn on techniques from various areas of algebraic geometry, see e.g., \cite{CG,IskovskikhManin,Pukhlikov,deF,kollarrat}. We now give a brief overview (for simplicity, over $\bC$) of progress on this problem over the past decade (see, for example, \cite{surveyschr} for a recent survey): these results were obtained using the specialization method introduced by Voisin,  developed by Colliot-Th\'el\`ene--Pirutka, and generalized by Schreieder, Totaro and others \cite{Voisinquartic,ColliotThelenePirutka,TotaroHypersurfaces, Schduke,  SchreiederSmallSlopes, SchreiederTorsion}.

Totaro \cite{TotaroHypersurfaces} applied the specialization method to show that a very general hypersurface of dimension $\dim(X)\geq 3$ and of degree at least $2\lceil \frac{\dim(X)+2}{3}\rceil$ is not stably rational over $\bC$.
In  \cite{SchreiederSmallSlopes}, Schreieder  improved this bound to $\log_2(\dim(X))+2$, based on a powerful version of the specialization method developed in \cite{Schduke}. 
Further improvements of this bound have been later obtained by Moe \cite{moe}, Nicaise--Ottem \cite{Nicaseottem},  Lange--Schreieder \cite{LSCH}, and, very recently, by Hotchkiss \cite{HotchkissDegreeFive} who showed the stable irrationality of a very general hypersurface of degree $\geq 5$ in every positive dimension.

In the case of quartic hypersurfaces, one has the following results. In dimension 3, stable irrationality of a very general quartic hypersurface was first proved by Colliot-Th\'el\`ene--Pirutka \cite{ColliotThelenePirutka}. Dimension 4 was covered by Totaro \cite{TotaroHypersurfaces}. Dimension 5 was settled by Nicaise--Ottem \cite{Nicaseottem} in characteristic 0 and  by Pavi\'c--Schreieder \cite{quart5sch} in characteristic different from 2. Recently, Ottem \cite{OttemQuarticSixfolds} settled dimension 6 in characteristic 0; his  argument is valid in any characteristic different from 2.

In this paper we settle the stable irrationality of very general quartic hypersurfaces in every positive dimension.

\begin{theo}[cf. Theorem~\ref{thm:mainreapt} and Corollary~\ref{coro:quartic}]\label{thm:main}
Let $k$ be an uncountable algebraically closed field with $\mathrm{char}(k)\ne 2$. Then a very general quartic hypersurface of any positive dimension over $k$ does not admit an integral decomposition of the diagonal, and in particular it is not stably rational.  
\end{theo}
Our proof  relies crucially on the  framework developed by Schreieder 
\cite{Schduke,SchreiederSmallSlopes} for the rationality problem of hypersurfaces.
Our result also generalizes to very general hypersurfaces of degrees $\geq5$, recovering the result of Hotchkiss   \cite{HotchkissDegreeFive}, see Remark~\ref{rema:higher}.

Theorem~\ref{thm:main} reduces the rationality problem of (very general) hypersurfaces to the case of cubics. It is a classical result of Clebsch that smooth cubic surfaces over an algebraically closed field are rational. Engel--de Gaay Fortman--Schreieder recently showed that very general complex cubic threefolds are not stably rational \cite{cubicmatroid}. The work of Katzarkov--Kontsevich--Pantev--Yu \cite{KKPY2025} is devoted to irrationality of very general cubic fourfolds. Very general cubics are believed to be stably irrational in all dimensions $\geq 4$. Note that the expectation concerns very general cubics, whereas \cite{TZc} supplies special stably rational smooth cubics.

 \subsection*{Overview of the proof}

Our proof of Theorem~\ref{thm:main} also applies the specialization method.  This method has two  steps. First, one needs to construct a {\it reference} variety, often  given by a (generic) quadric bundle carrying a certain unramified cohomology class. Second, one needs to realize a reference variety in a family  to deduce obstructions to stable rationality of the geometric generic fiber. 
In this paper we use the specialization method in the set up developed by
 Schreieder  for rationality problems of hypersurfaces \cite{Schduke,SchreiederSmallSlopes,SchreiederTorsion}. Schreieder's work   is fundamental for our application since we are using  families with a highly singular special fiber. We recall it in Section~\ref{sect:spec}.

 Here we briefly outline our construction. In Section 3, we construct  the quadric bundle data and the reference variety, which is closely related to the one in the work of Hotchkiss \cite{HotchkissDegreeFive}. The idea of the construction goes back to the method of Colliot-Th\'el\`ene and Ojanguren \cite{CTO} which has also been extensively used recently (in particular in \cite{Schduke}). We aim at constructing a fibration $f:Y\to S$  whose generic fiber is a quadric $Q$, and such that one has a nonzero unramified class in $\rH^m_{\mathrm{nr}}(k(Y)/k,\bZ/2)$ which satisfies certain properties (see Theorem \ref{prop:specialization}). In our case, the construction  adapts the strategy of  \cite{SchreiederSmallSlopes} which was inspired by \cite[Example 8]{HPT}.
 The quadric $Q$ is given by a diagonal form $\sigma \perp \langle F\rangle$ where $\sigma$ is a subform of a Pfister form $\Pi$ (see Proposition \ref{prop:quadric}), and $F$ is chosen in a particular way so that the unramified class can be constructed: in short  $\sigma \perp \langle F\rangle$ is not a subform of $\Pi$, but becomes a subform over some completions of $k(S)$.

In order to apply the specialization method and realize $Y$ in a family with generic fiber given by a hypersurface of small degree, one needs to    control   degrees of coefficients of $\sigma$, and the degree of $F$. Often $S$ is a rational variety with the function field $k(S)=k(z_1,\ldots z_m)$, and the coefficients of $\sigma$ are monomials in the $z_i$ or products of linear terms. This observation motivates the use of Segre embeddings in the next section.

In Section~\ref{sect:quart}, we explain the construction of the desired family of quartics. The main idea is to spread out the reference variety into a family, which is birational to a family of quartics. To achieve this, we first use the Segre embedding of a Segre variety $S$ (in our case, $S=(\bP^1)^a\times(\bP^2)^b$) 
to reduce the degrees of the coefficients of the generic quadric, which allows us to spread it out to (affine) polynomials of lower degrees. Then we separate coordinates on $S$ and on normal directions of $S$ to build a birational quartic family. The key step is to assign different weights on these coordinates and construct a quartic polynomial such that when we perform a weighted blowup, the leading terms of the polynomials become the generic quadric with coordinates in the normal directions. 

\medskip

\textbf{AI disclosure.} 
This project resulted from roughly two weeks of interaction mainly with ChatGPT-6 Astra. These interactions included questions related to another project (developed initially without use of AI), suggesting possible strategies, understanding and simplifying the proof. Claude and several models from the SAIR Open Math Model were also used during this process.

 The construction of a degeneration of quartics in dimension 10 using the Segre embedding of $\bP^2\times\bP^2$ was first found by ChatGPT.  The authors then guided ChatGPT to generalize this construction to increasingly higher dimensions (up to dimension 1015), before a general construction was  found in every dimension $\geq$ 7. 
 We wish to emphasize that the final constructions presented in Sections~\ref{s:Quadric} and~\ref{sect:quart} were found by ChatGPT, but the theoretical framework has been developed by humans. 
 The application of these constructions to stable irrationality rests on Schreieder’s work \cite{Schduke, SchreiederSmallSlopes, SchreiederTorsion}, which build on the specialization method of Voisin and Colliot-Thélène–Pirutka \cite{Voisinquartic,ColliotThelenePirutka}.

The paper was completely written by the authors, and the authors are responsible for its correctness. ChatGPT and Claude have been used for proofreading and finding literature in this process.

\medskip

\textbf{Acknowledgments. }Christian B\"ohning would like to thank the Simons Foundation for the inspiring research environment and excellent working conditions during a current stay in November and December of 2026. We would like to thank Stefan Schreieder and Yuri Tschinkel for many useful suggestions. Zhijia Zhang thanks James Hotchkiss for helpful conversations and explaining his work \cite{HotchkissDegreeFive}.

\section{Generalities}
\label{sect:spec}

For a field $k$ of characteristic different from $2$ we write $\rH^m(k,\bZ/2)$ for the Galois cohomology group.
For a field  extension $L/k$  we
write   $\rH^m_{\mathrm{nr}} (L/k,\bZ/2)$ for the unramified cohomology group defined as
$$\rH^m_{\mathrm{nr}} (L/k,\bZ/2)=\bigcap_v\mathrm{Ker}[\rH^m(L, \bZ/2)\stackrel{\partial^m_{v}}{\to}\rH^{m-1}(\kappa(v), \bZ/2)],$$ where the intersection is over all discrete valuations $v$ on  $L$ (of rank one), trivial on the field $k$, and where $\partial_{v}^m$ are residue maps (see \cite{CT95}).

Let $X$ be a proper integral variety defined over a field $k$. Recall that  $X$ has an {\it integral decomposition of the diagonal}   if there exist a
zero-cycle $z\in \CH_0(X)$ of degree $1$ and a cycle $\Gamma\in \CH_{\mathrm{dim}X}(X\times X)$ with support not dominating the first factor such that
$$\Delta_X=[X\times z]+\Gamma\in \CH_{\mathrm{dim}X}(X\times X).$$
Similarly, we say that an integral multiple $e\Delta_X$ for $e>0$ has an integral decomposition if the equality above holds with $e\Delta_X$ in place of $\Delta_X.$ 

The main degeneration theorem we use comes from \cite{Schduke, SchreiederSmallSlopes,SchreiederTorsion}.
The theorem below is a reformulation of \cite[Proposition 3.1]{SchreiederSmallSlopes} and \cite[Proposition 6.1]{SchreiederTorsion}, used in \cite[Theorem 2.2]{HotchkissDegreeFive}. Here, we state it in the form suitable for our application and we give a proof for convenience. 

\begin{theo}[Schreieder]\label{prop:specialization} 
Let $R$ be a discrete valuation ring with field of fractions $F$, and with an algebraically closed residue field $k$ of characteristic different from 2. Let $\mathcal X\to\Spec(R)$ be a proper and flat $R$-scheme
with integral geometric generic fiber. Let $\mathcal X_0$ be the special fiber and $Y$ an irreducible component of the reduction of $\mathcal X_0$, and assume that $\cX_0$ is reduced at the generic point of $Y$. Assume that there is a surjective morphism 
$$
f\colon Y\to S
$$
whose generic fiber is a smooth projective quadric $Q$ over $K=k(S)$, where $S$ is a proper integral $k$-variety  of dimension $m\ge1$. 
Assume that there exists a smooth integral closed subvariety $D\subset Q$, and a class $\alpha\in \rH^m(K, \bZ/2)$  such that 
\begin{enumerate}
\item the pullback $f^*(\alpha)$ to $k(Y)$ is nonzero and unramified: 
$$f^*(\alpha)\in \rH^m_{\mathrm{nr}}(k(Y)/k,\bZ/2);$$
\item the restriction  of $\alpha$ to $K(D)$ vanishes:
$$\alpha_{K(D)}=0\mbox{ in }\rH^m(K(D),\bZ/2);$$
\item the intersection of $Y$ and the singular locus of $\cX_0$ satisfies: $$(Y\cap\Sing(\cX_0))_{k(S)}\subset D.$$
\end{enumerate}
Then no odd multiple of the diagonal of the geometric generic fiber of $\cX$ admits
an integral decomposition, and in particular it is  not stably rational over $\overline{F}$.
\end{theo}

\begin{proof}
Let $\pi: Y'\to Y$ be an alteration of degree coprime to 2, with $Y'$ smooth,  which exists by results of de Jong and Gabber, see \cite{dejongaleration,gabber}.
By \cite[Proposition 6.1]{SchreiederTorsion}, it suffices to show that 
$$
\pi^*f^*(\alpha)\vert_{k(E)}=0
\text{ for every subvariety $E\subset \pi^{-1}( Y\cap\Sing(\cX_0))$}.
$$
Let $E\subset \pi^{-1}( Y\cap\Sing(\cX_0))$. If $\pi(E)$ does not dominate $S$ via $f$, then $\pi^*f^*(\alpha)\vert_{k(E)}=0$ follows from \cite[Theorem 9.2]{SchreiederSmallSlopes}. We assume that $\pi(E)$ dominates $S$. Then the generic point of $E$ over $k(S)$ maps to a scheme point $x\in D$ via $\pi$ by Assumption (3). By Assumption (2) and \cite[Lemma 2.3]{SchreiederTorsion}, we know that 
$$
\alpha_{\kappa(x)}=0\in \rH^m(\kappa(x),\bZ/2)
$$
where 
$\kappa(x)$ is the residue field  of $x$ and $\alpha_{\kappa(x)}$ is the pullback of $\alpha$ along $k(S)\hookrightarrow\kappa(x)$. Then $\pi^*f^*(\alpha)\vert_{k(E)}=0$ follows from the functoriality of Galois cohomology and field  inclusions 
$
k(S)\hookrightarrow \kappa(x)\hookrightarrow k(E).
$
\end{proof}
Such a triple $(\alpha,Q,D)$ as in Theorem~\ref{prop:specialization} is called a {\em Schreieder triple } in \cite{HotchkissDegreeFive}. 
In Section~\ref{sect:quart}, we apply Theorem~\ref{prop:specialization} to a family $\cX$ whose generic fiber is singular, but birational to the generic fiber of a degeneration of quartic hypersurfaces. 

We will also use the following lemma due to  Schreieder \cite[Lemma 2.4]{SchreiederSmallSlopes}.

\begin{lemm}\label{lemm:singdod}
    Let $k$ be a field. Let $X$ be a proper variety over $k$ and let $Y$ be a smooth proper variety  over $k$. Assume that $Y$ admits an integral decomposition of the diagonal and that $X$ is birational to $Y$ over $k$. Then $X$ also admits an integral decomposition of the diagonal. 
\end{lemm}

\section{A quadric fibration with a nontrivial unramified class}\label{s:Quadric}

In this section, we construct a quadric bundle satisfying conditions (1) and (2)  in Theorem~\ref{prop:specialization}. This will be a component of the special fiber of the degeneration family in Section~\ref{sect:quart}.

\subsection{Construction}

Let $k$ be an algebraically closed field of
characteristic $\mathrm{char}(k)\neq 2$.
Consider the Segre product variety $S=(\bP^1)^a\times(\bP^2)^b$ over
$k$, for $a\ge0$, $b\ge2$. We write its function field as
$$
 k(S)=k(u_1,\ldots,u_a,x_1,y_1,\ldots,x_b,y_b)
      =k(z_1,\ldots,z_m),\quad m=a+2b
$$
where  $u_i$ are affine coordinates on the $i$th $\bP^1$-factor
and $x_i,y_i$ are affine coordinates on the $i$th $\bP^2$-factor.
Throughout, we use the lexicographic order on monomials in
$z_1,\ldots,z_m$ given by
$$
 u_1>u_2>\cdots>u_a>x_1>y_1>x_2>y_2>\cdots>x_b>y_b.
$$

Let  
$$
\alpha=(z_1,\ldots,z_m)\in \rH^m(k(S),\mathbb Z/2).
$$ 

We  now explain the construction of a fibration in quadrics over $S$ that we will use as a {\it reference variety}. The generic fiber $Q$ will be constructed using subforms of  the Pfister form\footnote{We will use the notation $\Pf{z_1}=\langle 1, z_1 \rangle$ since $-1$ is a square in $k$.}
\begin{equation}\label{defpi}
\Pi=\Pf{z_1,\ldots,z_m}=\bigotimes_{i=1}^m\langle1,z_i\rangle,
\end{equation}
and one additional diagonal term.

Recall that $\Pi$ is a diagonal quadratic form of rank $2^m$ over $k(S)$. Let 
 $\Lambda\subset\Pi$ be the subform given by
\[
 \Lambda=\bigotimes_{i=1}^a\langle1,u_i\rangle
             \otimes\bigotimes_{j=1}^b\langle1,x_j,y_j\rangle.
\]
Put
\begin{equation}\label{eqn:P}
 P=\prod_{i=1}^a u_i\prod_{j=1}^b x_j
\end{equation}
and 
\[
 X=\prod_{j=2}^b x_j,\qquad
 Y=\prod_{j=2}^b y_j,\qquad D=x_1Y+y_1X.
\]
Let 
\begin{itemize}
\item [$\bullet$] $\cQ$ be the set of monomials in $z_1,\ldots,z_m$ which appear
as a coefficient of $\Lambda$ (including the constant $1$); 
\item[$\bullet$] $\cP$ be 
the subset of $\cQ$ obtained by deleting the following monomials
\[
 P,\quad\frac{P}{u_i},\quad\frac{Py_j}{x_j},\quad\frac{P}{x_j},
 \qquad i=1,\ldots,a,\quad j=1,\ldots,b;
\]
\item  $\cP'$ be the set of coefficients in $\Pi$ which are not in
$\cP\cup\{x_1y_1,XY\}$.
\end{itemize}
One can check that 
$$
|\cP|=2^a3^b-1-a-2b,\quad |\cP'|=2^m-|\cP|-2.
$$
We denote the $i$th element in $\cP'$ by $\mu_i$, under the increasing lexicographic
order. The order of $\mu_i$ plays no role in the construction. It is just a way to index elements in $\cP'$ for notational purposes.

\medskip

We write $\perp$ for the orthogonal sum of forms, e.g.,
$$\langle \mu_1\rangle\perp\langle \mu_2\rangle=\langle \mu_1,\mu_2\rangle .$$
Let
\[
 \sigma_0=\mathop{\perp}_{\mu\in\cP}\langle \mu\rangle.
\]
 For every integer $|\cP|\le r\le2^m-2$, we define
\begin{align}\label{eqn:sigmaprime}
     \sigma:=
 \displaystyle\mathop{\perp}_{i=1}^{r-|\cP|}\langle \mu_i\rangle
      \end{align}
Let 
\begin{equation}\label{eqn:sigmar}
\sigma_r:=\sigma\perp\sigma_0
\end{equation} be the rank $r$ quadratic form.

\begin{prop}\label{prop:quadric}
Let $a\geq 0, b\geq 2, m=a+2b$ be integers.
Let $$K_0=k(S)=k(z_1,\ldots,z_m),$$ and let $t$ be transcendental over
$K_0$. Put
$$
 \ell=\overline{k(t)},\qquad K=\ell(S)=\ell(z_1,\ldots,z_m).
$$
Fix an integer $r$ such that
$$2^a3^b-1-a-2b\le r\le2^m-2,$$ and let $\sigma_r$ be the 
quadratic form defined in~\eqref{eqn:sigmar}.
Let $f_1,f_2\in K_0^\times$ and put
$$
 F=f_1-t^2f_2\in K^\times.
$$
Assume the following conditions.
\begin{enumerate}
\item[(i)] There is an isometric embedding over $K_0$
\begin{equation}\label{eq:orthogonal-realization}
 \sigma_r\perp\langle f_1\rangle\subseteq\Pi.
\end{equation}
\item[(ii)] The form $\sigma_r$ represents $f_2$ over $K_0$.
\item[(iii)] For each discrete valuation $v$ on
 $K$, trivial on $\ell$, if $\Pi$ is
anisotropic over the completion $K_v$ of $K$ at $v$, and if $v(z_i)\ne0$ for some $i$, 
then
\[
 Ff_1\in K_v^{\times2}.
\]
\end{enumerate}
Let
\[
 \alpha=(z_1,\ldots,z_m)\in \rH^m(K, \mathbb Z/2),
\]
and let
\[
 Q=\{\sigma_r(U)+FV^2=0\}\subset\bP^r_K
\]
be a smooth projective quadric in coordinates $[U_1:\cdots:U_r:V]$ on $\bP^r_K$. 
 Let $\qquad H=\{V=0\}\subset\bP^r_K,$ and write $Q_H=Q\cap H$.
Then $Q_H$ is smooth and integral, and the following statements hold.
\begin{enumerate}
\item The pullback $\alpha_{K(Q)}$ is nonzero and belongs to
$\rH^m_{\nr}(K(Q)/\ell, \mathbb Z/2)$.
\item For every scheme  point $x$ of $Q_H$, the restriction
$\alpha|_{\kappa(x)}$ is zero. In particular, the specialization
of $\alpha_{K(Q)}$ at the generic point of $Q_H$ is zero.
\end{enumerate}
\end{prop}

\begin{rema}\label{rem:Ff1}
 The condition
\[
 Ff_1\in K_v^{\times2}
\]
is equivalent to $F/f_1=1-t^2f_2/f_1\in  K_v^{\times2}$.
A sufficient condition is
\[
 v(f_2/f_1)>0
\]
using Hensel's lemma and  $v(t)=0$.
\end{rema}

\begin{proof}[Proof of Proposition~\ref{prop:quadric}]

By \cite[Proposition 2.4]{HotchkissDegreeFive}, $\alpha_{K(Q)}$ is unramified.
We now show that $\alpha_{K(Q)}$ is nonzero.  The proof is similar to that of \cite[Proposition 6.1]{SchreiederSmallSlopes}. It is enough to establish that 
$$
\alpha_{K_E(Q_E)}\ne 0, 
$$ for every finite extension $E\subset\ell$ of $k(t)$, where $K_E=E(z_1,\ldots,z_m)$ and $Q_E$ is the quadric over $K_E$ defined by
$$
Q_E: \sigma_r(U)+(f_1-t^2f_2)V^2=0;
$$
the base change of $Q_E$ to $K$ is $Q$. Since $k$ is algebraically closed, a place of $E$ over $(t^{-1})$ extends to a discrete valuation $v$ on $E(z_1,\ldots,z_m)$ such that its residue field is $\kappa(v)=K_0=k(z_1,\ldots , z_m)$ and such that $ t^{-1}\in\mathfrak m$ is in the maximal ideal of the valuation ring $A$.

After the projective coordinate change $V'=tV$, the equation of $Q_E$ is
$$
 \sigma_r(U)+(t^{-2}f_1-f_2)(V')^2=0.
$$
Every coefficient of $\sigma_r$ and the element $f_2$ are nonzero
elements of $K_0$, hence units of $A$. It follows that the last coefficient $(t^{-2}f_1-f_2)$ is also
a unit, with reduction $-f_2$ in the residue field $K_0$. The same equation therefore defines
a smooth projective quadric $\mathfrak Q$ over $A$.

The  special fiber of  $\mathfrak Q$ 
is the smooth geometrically integral quadric
\[
 Q_\infty=\{\sigma_r(U)-f_2(V')^2=0\}\subset\bP^r_{K_0}.
\]
Hence the local ring of $\mathfrak Q$ at the generic point of $Q_{\infty}$ is a discrete valuation ring on $E(z_1,\ldots , z_m)(Q)$ with residue field $\kappa(w)=K_0(Q_\infty)$. We denote by $w$ the corresponding valuation.
Hypothesis~(ii) gives a $K_0$-point of $Q_{\infty}$ 
with $V'=1$,  hence  $
 K_0(Q_\infty)$ is a purely transcendental extension of $K_0=k(z_1,\ldots , z_m)$.
The restriction $\alpha|_{\kappa(w)}$ is then the class of  $(z_1,\ldots,z_m)$ in $K_0(Q_\infty)$, which is nonzero. Hence $\alpha_{K_E(Q_E)}$ is nonzero as claimed.

\smallskip
We now establish part (2) of the Proposition.
Let $x$ be any scheme point of $Q_H$. A nonzero representative
of its homogeneous coordinates is an isotropic vector of $\sigma_r$
over $\kappa(x)$. The embedding
\eqref{eq:orthogonal-realization} makes $\Pi$ isotropic over this
field, hence
\[
 \alpha|_{\kappa(x)}=0
\]
as claimed.
\end{proof}

\subsection{Concrete choices for $f_1$ and $f_2$}

Such choices of $f_1$ and $f_2$ satisfying the conditions in Proposition~\ref{prop:quadric} are generally not unique. Here we present a concrete choice that we will use in Section~\ref{sect:quart}. 
\begin{prop}[Explicit coefficients]\label{prop:explicit-coefficients}
With the notation of Proposition~\ref{prop:quadric}, define polynomials in the affine
coordinates of $S$ by
\[
 G=\prod_{i=1}^a(1+u_i)\prod_{j=1}^b(1+x_j+y_j)
\]
and
\begin{equation}\label{eqn:f1f2}
 f_1=G^2x_1XD=G^2x_1X(x_1Y+y_1X),\qquad
 f_2=x_1^2y_1X^2Y\prod_{i=1}^a u_i.
\end{equation}
These polynomials satisfy hypotheses~(i)--(iii) of
Proposition~\ref{prop:quadric} for every integer $r$ with 
$2^a3^b-1-a-2b\le r\le2^m-2$, $a\geq 0,b\geq 2.$

Moreover, $f_1$ has degree two in each $u_i$ and total degree four
in each pair $(x_j,y_j)$, while $f_2$ has degree one in each $u_i$
and total degree three in each pair $(x_j,y_j)$.

\end{prop}

\begin{rema}
By the computation of degrees,  $f_1-t^2f_2$ has degree at most four in each factor of $S$; we will see that it admits an ambient quartic lift  constructed in Section~4.
\end{rema}

\smallskip
\begin{exam}
    For $a=1$ and $b=2$, the formulas give
\[
\begin{aligned}
 X&=x_2,\qquad Y=y_2,\qquad D=x_1y_2+y_1x_2,\\
 G&=(1+u_1)(1+x_1+y_1)(1+x_2+y_2),\\
 f_1&=(1+u_1)^2(1+x_1+y_1)^2(1+x_2+y_2)^2
             x_1x_2(x_1y_2+y_1x_2),\\
 f_2&=u_1x_1^2y_1x_2^2y_2.
\end{aligned}
\]
The total degrees of $f_1$ and $f_2$ are $10$ and $7$, respectively.
Their degrees in $u_1$ are $2$ and $1$, and their total degrees in
each pair $(x_j,y_j)$ are $4$ and $3$.
\end{exam}

\begin{proof}[Proof of Proposition~\ref{prop:explicit-coefficients}]
We recall that, by definition, the distinct squarefree
monomials $x_1y_1$ and $XY$ both lie outside $\cQ$.

We choose a basis for the underlying vector space of $\Pi$ as follows: for each squarefree monomial $h$ in the variables $z_i$, denote by
$e_h$ the corresponding vector of $\Pi$. In this basis, $\Pi$ has diagonal form as in  (\ref{defpi}). 
 Define
vectors
\[
 \xi=Y e_{x_1y_1}-y_1e_{XY},
 \qquad
 \rho=G\bigl(Xe_{x_1y_1}+x_1e_{XY}\bigr).
\]
They are orthogonal, and one has
\[
\begin{aligned}
 \Pi(\xi)&=x_1y_1Y^2+XYy_1^2=y_1YD,\\
 \Pi(\rho)&=G^2(x_1y_1X^2+XYx_1^2)=f_1.
\end{aligned}
\]
Every basis vector indexed by $\cP\cup\cP'$ is orthogonal to both
vectors. Mapping the diagonal basis vectors of $\sigma_r$ to the
corresponding $e_h$
hence gives 
\[
 \sigma_r\perp\langle f_1\rangle\ \subseteq\ \Pi,
\]
with the last summand represented by $\rho$.

The monomial
\[
 h=y_1Y\prod_{i=1}^a u_i
\]
belongs to $\cP$. Indeed, it differs from $P$ in all $b\ge2$
$\bP^2$-factors.  Hence  $\sigma_r$ represents
\begin{equation}\label{eq:f2-represented}
 f_2=(x_1X)^2h=\sigma_r(x_1Xe_h).
\end{equation}
Hence (i) and (ii) are satisfied, and (iii) follows from  
Lemma~\ref{lem:valuation-estimate} below and Remark \ref{rem:Ff1}.

For the degrees, observe that $G$ has degree one in each $u_i$ and
total degree one in each pair $(x_j,y_j)$.
The factor $x_1XD$ has total degree two in each pair $(x_j,y_j)$
and is independent of the $u_i$. The assertion about the degrees of $f_1$ follows.
For $f_2$ it is enough to observe
\[
 f_2=\left(\prod_i u_i\right)\left(\prod_j x_j^2y_j\right).
\]

\end{proof}

\begin{lemm}
\label{lem:valuation-estimate}
Let $f_1,f_2,G$ be the polynomials of
Proposition~\ref{prop:explicit-coefficients}. Let $v$ be a discrete valuation on $K=\ell(S)$,  trivial on $\ell$,
and suppose that $\Pi$ is anisotropic over $K_v$. 
 If at least one of the base coordinates
$u_i,x_j,y_j$ has nonzero valuation, then $v(f_2/f_1)>0$. 
\end{lemm}

\begin{proof}
Recall that $D=x_1Y+y_1X$. We claim that
\[
\begin{aligned}
 v(1+u_i)&=\min\{0,v(u_i)\},\\
 v(1+x_j+y_j)&=\min\{0,v(x_j),v(y_j)\},\\
 v(D)&=\min\{v(x_1Y),v(y_1X)\}.
\end{aligned}
\]
Indeed, 
the diagonal forms
$
 \langle 1,u_i\rangle, \langle 1,x_j,y_j\rangle,
 $
 and
 $
 \langle x_1Y,y_1X\rangle
$
are subforms of $\Pi$, hence anisotropic over $K_v$. For any of these forms,
denote its coefficients by $a_1,\ldots,a_s$, ordered so that
\[
 v(a_1)=\min_i v(a_i).
\]
If $v(a_1+\cdots+a_s)>v(a_1)$, applying Hensel's lemma to the polynomial
\[
 h(T)=T^2-1+\frac{a_1+a_2+\cdots+a_s}{a_1}\in\mathcal O_v[T],
\]
we get $c\in K_v^\times$ such that
$
 a_1c^2+a_2+\cdots+a_s=0,
$
contradicting that $\Pi$ is anisotropic.

Put
\[
\begin{gathered}
 c_i=v(u_i),\qquad p_j=v(x_j),\qquad q_j=v(y_j),\\
 A=v(x_1Y),\qquad B=v(y_1X).
\end{gathered}
\]
We deduce from the above that
\[
 v(G)=\sum_i\min(0,c_i)+\sum_j\min(0,p_j,q_j).
\]
Since
\[
 \frac{f_2}{f_1}
 =\frac{\bigl(\prod_i u_i\bigr)
          \bigl(\prod_j x_jy_j\bigr)}{G^2D},
\]
we have that
\begin{multline*}
 v(f_2/f_1)
 =\sum_i c_i+\sum_j(p_j+q_j)-2\sum_i\min(0,c_i)
    -2\sum_j\min(0,p_j,q_j)-\\-\min(A,B)
     =\sum_i|c_i|+\max(A,B)-2\sum_j\min(0,p_j,q_j),
\end{multline*}
where the second equality uses the fact that $A+B=\sum_j(p_j+q_j)$. We deduce
 $
  v(f_2/f_1)\ge\sum_i|c_i|+
       \frac12\sum_j\bigl(|p_j|+|q_j|\bigr),
$
which is strictly positive by assumption.
\end{proof}

\section{Quartic degeneration}\label{sect:quart}

 In this section, we construct the desired degeneration of quartic hypersurfaces.  
Let $m=a+2b$, for some $a\geq 0, b\geq 2$, and 
$$
 r_0=2^a3^b-1-a-2b,\quad r_0\leq r\leq2^m-2,\quad n=r-r_0.
 $$
 Recall that $\cQ$ is the set of monomials in $z_1,\ldots,z_m$ which appear
as a coefficient in the quadratic form  $\Lambda=\bigotimes_{i=1}^a\langle1,u_i\rangle
             \otimes\bigotimes_{j=1}^b\langle1,x_j,y_j\rangle$ (including the constant $1$).
             Consider the standard Segre embedding 
$$
S=(\bP^1)^a\times(\bP^2)^b\hookrightarrow\bP^{2^a3^b-1},
$$
given by the homogenization of monomials in $\cQ$.
We label the corresponding coordinates on the target projective space by $W_q$, $q\in \cQ$. 

Further embed it into a larger projective space 
$$
\varphi: S\hookrightarrow\bP^{m+r},
$$
with coordinates $[(W_q)_{q\in\cQ}:s_1:\cdots:s_n]$. Note that $ |\mathcal{Q}| +n -1=m+r$.  The image $\varphi(S)$ is contained in the linear section $
\{s_1=\cdots=s_n=0\}$. 

The rest of this section is devoted to an explicit construction of the quartic family, which involves a lot of computations. Although some of the properties of the Segre embedding are known (e.g. in Sections \ref{section:liftlin} and \ref{sect:sigma_0}), we write all the constructions in detail.   We outline the main ideas first. 

\subsection{Strategy}

The main goal is to realize the reference variety $Y\to S$ with generic fiber the quadric $Q$ constructed in Proposition~\ref{prop:quadric} in a family whose generic fiber is birational to a smooth quartic in $\{F_\tau=0\}\subset\bP^{m+r}$. On an affine chart $\bA^{m}\times\bA^r$, we identify $\bA^m$ with an affine chart of $S$, and $\bA^r$ with its normal directions with respect to the embedding $\varphi$. To obtain a projective family in the end, we simply take the projective closure in $S\times \bP^r$. This is convenient for the computation although it results in a singular generic fiber. 

To make the quadric $Q$ appear in a component  of the central fiber, the main idea is to assign appropriate positive weights to coordinates in the normal directions, and zero weights to $\bA^m$, so that after a weighted rescaling of the coordinates (which can be seen as a chart of a weighted blowup), the central fiber is given by the least weighted terms of $F_\tau$.

For each term of  $Q$, we find a quartic polynomial on $\bP^{m+r}$ such that in an affine chart, its least weighted term is quadratic in the normal directions of $S$, and is exactly the term of $Q$. Recall $Q$ is given by a quadratic form $\sigma_0\perp\sigma \perp\langle f_1-t^2f_2\rangle$. For terms in $\sigma$ and $f_1-t^2f_2$, this becomes simple due to the Segre embedding, which reduces the degree of a monomial on $S$ to the maximal degree on a $\bP^1$ or $\bP^2$-factor (e.g., a degree 4 term $u_1x_1^2y_2$  becomes a quadratic term $W_{u_1x_1y_2}W_{x_1}$). 

By the construction in Section~\ref{s:Quadric}, $f_1-t^2f_2$ is quartic in each factor, and terms in $\sigma_0$ and $\sigma$ are at most quadratic in each factor. We will assign the linear variables $s_1,\ldots,s_n$ in the normal directions  to $\sigma$. Since the defining functions of the remaining normal coordinates are nonlinear in the original ambient coordinates, the main difficulty lies in finding quartic polynomials singular along $\varphi(S)$, and whose least weighted term on this affine chart is a term in $\sigma_0$. This is achieved in Sections~\ref{sect:sigma_0} and~\ref{sect:weight}, with the quartics chosen being products of quadratic elements in the ideal sheaf of $\varphi(S)$.

 \subsection{Lifting coefficients in linear variables}\label{section:liftlin}

In the affine chart of $\bP^{m+r}$ given by $\{W_1\ne0\}$, using coordinates $w_q=W_q/W_1$, one sees that 
\begin{align}\label{eqn:Scirc}
    S^\circ:=\varphi(S)\cap \{W_1\ne 0\}\simeq \bA^{m}
\end{align}
with coordinates $w_{u_i},w_{x_j},w_{y_j}$, since other coordinates are products of these. 
Let $\{\mu_i\}$ be the coefficients in $\sigma$ as is defined in~\eqref{eqn:sigmaprime}. Note that for each $\mu_i$, and for each factor ($\bP^1$ or $\bP^2$) of $S$, one has that $\mu_i$ is at most quadratic in variables $\{u_i\}$ or $\{x_j, y_j\}$ of this factor.  Hence we can find a homogeneous quadratic polynomial $Q_d$ in variables $(W_q)_{q\in \cQ}$ such that when restricting to the affine chart $ S^\circ$, we have that
\begin{align}
    \label{eqn:Qiconst}
    Q_d\vert_{S^\circ}(u_1,\ldots u_a,x_1,y_1,\ldots,x_b,y_b)=\mu_d.
\end{align}
This means that, after writing the restriction of $Q_d$ to $S^\circ$  in variables $w_{u_i}, w_{x_j}, w_{y_j}$, and evaluating it at 
$$
w_{u_i}\mapsto u_i,\quad w_{x_j}\mapsto x_j,\quad w_{y_j}\mapsto y_j, 
$$
we obtain the polynomial $\mu_d$.
 Similarly, by Proposition \ref{prop:explicit-coefficients}, we can find a homogeneous quartic polynomial $F$ in $(W_q)_{q\in \cQ}$ such that 
\begin{align}
    \label{eqn:Fconst}
F\vert_{S^\circ}(u_1,\ldots u_a,x_1,y_1,\ldots,x_b,y_b)=f_1-t^2f_2
\end{align}
where $f_1,f_2$ are defined as in~\eqref{eqn:f1f2}.

\subsection{Lifting coefficients of $\sigma_0$}\label{sect:sigma_0}
In this subsection, we construct a quartic polynomial, for each $\mu\in\cP$, that later recovers, in a way that will be made more precise below, the term in $\sigma_0$ with coefficient $\mu$ in the central fiber of the degeneration.

Let $\mu\in\cP$. Note that $\mu$ is a product 
\begin{align}
    \label{eqn:beta}
    \mu=\beta_1\cdots\beta_a\cdot \beta_{a+1}\cdots\beta_{a+b},
\end{align}
where $\beta_i\in\{1,u_i\},\beta_{a+j}\in\{1,x_j,y_j\}$, for $i=1,\ldots,a$ and $j=1,\ldots,b$. Put
$$
i_0=i_0(\mu):=\max\{i:\beta_{i}=1 \text{ or } i>a \text{ and }\beta_i=y_{i-a} \}.
$$
By the definition of $\cP$, every element in $\cP$ differs from $P$ defined in~\eqref{eqn:P} by at least 2 elements, which implies that $i_0\geq2$. We create a matrix $M=M^{(i_0)}$ such that 
\begin{itemize}
    \item it has two or three columns which are labeled (in an increasing order) by $\{1,u_{i_0}\}$ if $i_0\leq a$, and by $\{1,y_{i_0-a},x_{i_0-a}\}$ if $i_0>a$ (i.e., affine coordinates including the constant 1 coming from the $i_0$th factor of $S$);
    \item its rows are labeled (in an increasing order) by all other elements in $\cQ$ which have no factors from the affine coordinates on $i_0$th factor of $S$ (i.e., elements in $\cQ$ not containing any $u_{i_0}$ if $i_0\leq a$, or any $x_{i_0-a}, y_{i_0-a}$ if $i_0>a$);
    \item the $(i,j)$th entry $M_{i,j}$ of $M$ is given by $W_\nu$ where $\nu$ is the product of the labels of the $i$th row and $j$th column.
\end{itemize}
The entries of $M$ are simply all of $(W_{q})_{q\in \cQ}$, but organized in a way depending on $i_0$. We note that the matrix $M^{(i_0)}$ is called a {\em tensor flattening}, see e.g., \cite[Chapter 9.2]{sturmfels}

\begin{exam}
    Assume that $S=(\bP^1)^2\times\bP^2$ (we use this only as a low-dimensional example to illustrate the construction of $M^{(i_0)}$; in our hypothesis $b\geq 2$). Then we have
    $$
    M^{(2)}=\begin{pmatrix}
        W_1&W_{u_2}\\
        W_{y_1}&W_{u_2y_1}\\
        W_{x_1}&W_{u_2x_1}\\
        W_{u_1}&W_{u_1u_2}\\
        W_{u_1y_1}&W_{u_1u_2y_1}\\
        W_{u_1x_1}&W_{u_1u_2x_1}
    \end{pmatrix},\quad M^{(3)}=\begin{pmatrix}
        W_1&W_{y_1}&W_{x_1}\\
        W_{u_2}&W_{u_2y_1}&W_{u_2x_1}\\
        W_{u_1}&W_{u_1y_1}&W_{u_1x_1}\\
        W_{u_1u_2}&W_{u_1u_2y_1}&W_{u_1u_2x_1}
    \end{pmatrix}.
    $$
\end{exam}
Now we find the location of $W_\mu$ in $M=M^{(i_0)}$ and write
$$W_\mu=M_{i,j}.$$  
The quartic polynomial associated with $\mu$ is a product of two $2\times2$ minors, which below we refer to as $D_1$ and $D_2$:
\begin{align}
    \label{eqn:minors}
    F_\mu=
\det\begin{pmatrix}
W_1 & M_{1,j+1}\\
M_{i+1,1} & M_{i+1,j+1}
\end{pmatrix}
\det\begin{pmatrix}
W_\mu & M_{i,j+1}\\
M_{i+1,j} & M_{i+1,j+1}
\end{pmatrix}.
\end{align}
We observe that $F_\mu$ is well-defined, and $\{F_\mu=0\}$ is singular along $\varphi(S)$.
\begin{lemm}\label{lemm:welldefine}
    For each $\mu\in \cP$, the following holds:
    \begin{enumerate}
        \item The coordinate $W_{\mu}$ corresponding to the monomial $\mu$ is not in the last row or column in $M^{(i_0)}$. In particular $F_\mu$ is well-defined.
        \item We have $F_\mu\in H^0(\bP^{m+r},\cI_S^2 (4))$ where $\cI_S$ is the ideal sheaf of $\varphi(S)\subset \bP^{m+r}$.
    \end{enumerate} 
\end{lemm}
\begin{proof}
    (1) follows from the construction of $i_0$. Indeed, since $\beta_{i_0}$ is not divisible by $u_{i_0}$ or $x_{i_0-a}$, there is a column to the right of $\mu$. On the other hand, by the definition of $\cP$, every element in $\cP$ differs from $P$ defined in~\eqref{eqn:P} by at least 2 elements, which implies that there is an $i<i_0$ such that $\beta_i$ is not divisible by $u_{i}$ or $x_{i-a}$. Since the $i$th factor is larger in lexicographic order than the $i_0$th factor, we know that $\mu$ is not in the bottom row. 

    To show (2), it suffices to notice that $M^{(i_0)}$ has rank 1 when restricted to $\varphi(S)$ and thus every $2\times 2$ minor vanishes along $\varphi(S)$. Indeed, every row is a multiple of the first row on $S^\circ$.  
    \end{proof}

\subsection{Change of coordinates and weights.}
\label{sect:weight}
Let 
$$
\cN=\{q\in\cQ| \,\,q\text{ has degree $\geq2$ as a polynomial in $z_1,\ldots,z_m$}\}.
$$

We observe that there is a bijection 
\begin{align}\label{eqn:bij}
    g: \cP\to \cN
\end{align}
defined as follows: let $\mu\in \cP$ be given. Find $M_{i,j}=W_\mu$. We define $g(\mu):=\nu$ where $W_\nu=M_{i+1,j+1}$. It is well-defined by Lemma~\ref{lemm:welldefine} and its image lies in $\cN$ since every linear polynomial in $\cP$ lies on the first row or first column of $M$. 
   We can construct the inverse as follows. For any $\nu\in\cN$, we write 
    $$
    \nu=\nu_1\cdots\nu_a\nu_{a+1}\cdots\nu_{a+b}
    $$
    where  $\nu_i\in\{1,u_i\},\nu_{a+j}\in\{1,x_j,y_j\}$, for $i=1,\ldots,a$ and $j=1,\ldots,b$. Put 
    $$
    j_0=\max\{i:\nu_i\ne 1\}.
    $$
Find $
W_{\nu}=M^{(j_0)}_{i,j}$. One can check that $g^{-1}(\nu)=\mu$ where $W_\mu=M^{(j_0)}_{i-1,j-1}$.

Now we write down the change of variables realizing~\eqref{eqn:Scirc}. As before, we denote the standard affine coordinates on the chart $\{W_1\ne 0\}\cap\bP^{m+r}$ by  
$$
w_q:=\frac{W_q}{W_1},\quad q\in \cQ\setminus\{1\}\,\qquad s'_i:=\frac{s_i}{W_1}, \quad i=1,\ldots,n,$$ and set $w_1=1$.

We apply the change of variables 
\begin{equation}\label{eq:varch}
\eta_q:=w_q-\prod_{z_i \mid q}w_{z_i},\quad \text{ for all } q\in\cN,
\end{equation}
where the product runs through linear terms arising in $q$. For example, $\eta_{u_1x_1y_2}=w_{u_1x_1y_2}-w_{u_1}w_{x_1}w_{y_2}$.
Then in the coordinates 
\begin{align}
    \label{eqn:coord}
    (w_{u_i},w_{x_j},w_{y_j},\eta_q,s'_l) 
\end{align}
 of $\bA^{m+r}$,
for
$i\in\{1,\ldots,a\},j\in\{1,\ldots,b\},l\in\{1,\ldots,n\}, q\in\cN,
$ 
we have that 
\begin{align*}
    S^\circ=\{\eta_q=s_l'=0\}\simeq \bA^m.
\end{align*}

\medskip

To put the quartics $F_\mu$ into an appropriate family, we assign weights to the coordinates as follows. We identify $u_i=w_{u_i}, x_j=w_{x_j}, y_j=w_{y_j}$ on the affine chart and we define
$$
\mathrm{wt}(u_i)=\wt(x_j)=\wt(y_j)=0,\quad \mathrm{wt}(s_l')=1.
$$
We recall that $r_0=2^a3^b-a-2b-1$. Assume that elements in $\cN$ are in increasing lexicographic order $n_1<n_2<\ldots<n_{r_0}$. 
The weights of the
corresponding coordinates are given by
\[
 \wt(\eta_{n_i})=r_0+1-i.
\]

\begin{lemm}\label{lemm:weightFmu}
Let $g$ be the bijection defined  in~\eqref{eqn:bij}.
    For each $\mu\in\cP$, the restriction of $F_\mu$ to $\{W_1\ne 0\}\cap\bP^{m+r}=
    \bA^{m+r}$ in coordinates~\eqref{eqn:coord} has the form 
    $$
    F_\mu\vert_{\bA^{m+r}}=\mu\cdot \eta_{g(\mu)}^2+\text{terms with higher weights than $\eta^2_{g(\mu)}$}.
    $$
\end{lemm}
\begin{proof}
 In this proof, we abuse the notation $\eta_{q}$ and write
$$
\eta_q=\begin{cases}
    \eta_q&\text{if }q \in\cN \text{ i.e., }\deg(q)\geq2, \\
   0 &\text{if }q\not\in\cN \text{ i.e., }\deg(q)\leq1. 
\end{cases}
$$
    First, we look at the second minor in~\eqref{eqn:minors}
$$
D_2=\det\begin{pmatrix}
W_\mu & M_{i,j+1}\\
M_{i+1,j} & M_{i+1,j+1}
\end{pmatrix}=W_\mu W_{g(\mu)}-W_{\nu_2\mu_1} W_{\nu_1\mu_2}
$$
where
\begin{itemize}
    \item $\mu_1$, $\nu_1$ are the labels of the $i$th and the $(i+1)$th rows of $M$,
    \item  $\mu_2$, $\nu_2$ are the labels of the $j$th and the $(j+1)$th columns of $M$; note that $\mu_2$, $\nu_2$ are then at most linear polynomials in $u_i, x_j, y_j$.
\end{itemize}
Note that $\mu=\mu_1\mu_2$ and $g(\mu)=\nu_1\nu_2.$
Using \eqref{eq:varch}, we compute that 
$$
D_2\vert_{\bA^{m+r}}=\mu\eta_{g(\mu)}+g(\mu)\eta_\mu+\eta_{\mu}\eta_{g(\mu)}-\eta_{\nu_2\mu_1}\eta_{\nu_1\mu_2}-\nu_2\mu_1\eta_{\nu_1\mu_2}-\nu_1\mu_2\eta_{\nu_2\mu_1}.
$$
Since $W_{g(\mu)}$ is the bottom-right element of $D_2$, we know that 
\begin{align}\label{eqn:var}{g(\mu)}>\mu,\nu_1\mu_2,\nu_2\mu_1\quad\text{in lexicographic order}.
\end{align}
If all of them are in $\cN$, this implies that
$$
\wt(\eta_{g(\mu)})<\wt(\eta_\mu),\wt(\eta_{\nu_1\mu_2}),\wt(\eta_{\nu_2\mu_1}),
$$
and we conclude that the term in $D_2\vert_{\bA^{m+r}}$ with the least weight is $\mu\eta_{g(\mu)}$. If one of these monomials is not in $\cN$, the corresponding term in $D_2\vert_{\bA^{m+r}}$ vanishes and the same conclusion holds.

A similar computation shows that the term in $D_1\vert_{\bA^{m+r}}$ with the least weight is $\eta_{g(\mu)}$. Multiplying $D_1$ and $D_2$ completes the proof. 
\end{proof}
\begin{rema}
    We note that any quartic polynomial satisfying the conditions in Lemma~\ref{lemm:weightFmu} would work for our degeneration construction. The choices of $F_\mu$ are thus not unique and the one defined by~\eqref{eqn:minors} is only one of them.
\end{rema}

\subsection{Degeneration of quartics}

Finally, we are ready to construct the degeneration of quartics. Let $L=\overline{k(t)}$ and $R=L[[\tau]]$ be a discrete valuation ring. Using definitions \eqref{eqn:Qiconst} and \eqref{eqn:Fconst}, we view $F_\mu,  F$ as quartic forms, and $Q_j$ as quadratic forms on $\bP^{m+r}$ over $R$, with no variables $s_1,\ldots,s_n$. For each $\mu\in \mathcal P$, we define 
$$
\omega_\mu:=\mathrm{wt}(\eta_{g(\mu)}).
$$
Let $J$ be a {\em general} quartic form on $\bP^{m+r}$ over $R$, and define a family of quartics given by 
\begin{align}\label{eqn:Ftau}
F_\tau
= \sum_{\mu\in\mathcal{P}}
   \left(\tau^{2(r_0-\omega_\mu)}F_\mu\right)
   + \tau^{2(r_0-1)}\sum_{j=1}^{n}\left(s_j^2Q_j\right)
+ \tau^{2r_0}F
   + \tau^{2r_0+1}J.
\end{align}

The powers of $\tau$ are chosen so that after passing to affine coordinates~\eqref{eqn:coord}, and assigning weight $1$ to $\tau$, we have the following:
$$
\begin{tabular}{|c|c|}
\hline
    $f$ &  \text{least weight of a term in }$f$\\\hline
      $ \tau^{2(r_0-\omega_\mu)}F_\mu$&$2r_0$ \\\hline $\tau^{2(r_0-1)}s_j^2Q_j$&$2r_0$\\\hline 
$\tau^{2r_0}F$&$2r_0$\\\hline
$\tau^{2r_0+1}J$&$>2r_0$\\\hline
\end{tabular}
$$
The first row follows from Lemma~\ref{lemm:weightFmu}; the second row follows from~\eqref{eqn:Qiconst}; the third row follows from~\eqref{eqn:Fconst}; the last row is obvious.
Let 
\begin{align}
    \label{eqn:cZ}
    \cZ=\{F_\tau=0\}\subset\bP^{m+r}_R.
\end{align}
Then 
$$
\cZ\to \Spec(R)
$$ is a proper and flat morphism, whose generic fiber is a smooth quartic hypersurface in $\bP^{m+r}$ over $L(\!(\tau)\!)$. The central fiber $\cZ_0$ of $\cZ$ is a double quadric cone given by 
$$
\cZ_0=\{F_{g^{-1}(n_1)}=0\}\subset\bP^{m+r}_L,\quad n_1\text{ is the smallest element in }\cN.
$$
Now we construct a degeneration 
$$
\cX\to \Spec(R)
$$
such that the generic fibers of $\cX$ and $\cZ$ are birational over $L(\!(\tau)\!)$, and $\cX$ satisfies the conditions of Theorem~\ref{prop:specialization}.
We first take the generic fiber $\cZ_\eta$ of $\cZ$ and restrict it to the affine chart 
$$
\cZ_\eta^\circ=\cZ_\eta\cap\{W_1\ne 0\}\subset\bA^{m+r}_{L(\!(\tau)\!)}.
$$
Write $f_\tau$ for the dehomogenization of $F_\tau$ with respect to $W_1$. In the coordinates 
$$
 (w_{u_i},w_{x_j},w_{y_j},\eta_{g(\mu)},s'_l), 
$$
the polynomial $f_\tau$ takes the form 
$$
f_\tau=\sum_{\mu\in\mathcal{P}}
   \tau^{2(r_0-\omega_\mu)}\mu\eta_{g(\mu)}^2
   + \tau^{2(r_0-1)}\sum_{j=1}^{n}\mu_j\cdot(s'_j)^2
+ \tau^{2r_0}(f_1-t^2f_2)
   + h
$$
where $h$ is a polynomial each term of which has combined weight at least $2r_0+1$ ($\tau$ has weight 1).
On $\bA^{m+r}_{L(\!(\tau)\!)}$, we introduce new coordinates 
$$
w_{u_i},\quad w_{x_j},\quad w_{y_j},\quad T_{g(\mu)}:=\tau^{-\omega_\mu}\eta_{g(\mu)},\quad T_{l}:=\tau^{-1} s_l',
$$
for $\mu\in\cP, l\in\{1,\ldots,n\}.$ Note that this is simply scaling variables $\eta_{g(\mu)}$ and $s_l'$, so it gives an isomorphism on $\bA^{m+r}_{L(\!(\tau)\!)}$. Since each term in $h$ has weight at least $2r_0+1$, and $w_{u_i},w_{x_j}, w_{y_j}$ have weight $0$, we know that 
$$
h(w_{u_i},w_{x_j},w_{y_j},\tau^{\omega_\mu}T_{g(\mu)},\tau T_l)=\tau^{2r_0+1}h'(w_{u_i},w_{x_j},w_{y_j},T_{g(\mu)}, T_l)
$$
for some polynomial $h'$. 
In the new coordinates 
$$
(w_{u_i},w_{x_j},w_{y_j},T_{g(\mu)}, T_l),
$$
we have that 
$$
f_\tau=\tau^{2r_0}\left(\sum_{\mu\in\mathcal{P}}\mu T
   _{g(\mu)}^2
   + \sum_{j=1}^{n}\mu_j T_j^2
+ (f_1-t^2f_2)
   + \tau h'\right),
$$
and 
\begin{align}\label{eqn:Zcirc}
\cZ^\circ_\eta=\left\{\sum_{\mu\in\mathcal{P}}\mu T
   _{g(\mu)}^2
   + \sum_{j=1}^{n}\mu_jT_j^2
+ (f_1-t^2f_2)
   + \tau h'=0\right\}\subset\bA_{L(\!(\tau)\!)}^{m+r}.
\end{align}
Note that after substituting $\tau=0$ in~\eqref{eqn:Zcirc}, we obtain an affine chart of the quadric $Q$ given in Proposition~\ref{prop:quadric}. Recall that $S^\circ\simeq \bA^m$. We take a projective model $\cX_\eta$ of $\cZ_\eta^\circ$ by taking its projective closure in $S\times\bP^r$ under the natural embedding
$$
\bA^{m+r}\hookrightarrow S\times\bP^r
$$
given in coordinates by 
\begin{multline*}
    (w_{u_i},w_{x_j},w_{y_j},T_{g(\mu)}, T_l)\mapsto \left(\prod_{i=1}^a(1:u_i)\times\prod_{j=1}^b(1:x_j:y_j)\right)\times(1:T_{g(\mu)}:T_l).
\end{multline*}
Let 
$$
\cX\subset S\times\bP^r
$$
be the corresponding projective model over $\Spec(R)$ given by the same equation.

\begin{rema}
    We note that $\cX_\eta$ is   singular but geometrically integral. The scaling of affine coordinates above can also be seen as a chart of a weighted blowup, which provides a way to construct a proper  model for our degeneration with smooth generic fiber. For simplicity, we decide not to do so and proceed with taking the projective closure, since Theorem~\ref{prop:specialization} does not require smoothness of the generic fiber.
\end{rema}
\begin{prop}\label{prop:degen}
    The morphism $\cX\to\Spec(R)$ satisfies the assumptions in Theorem~\ref{prop:specialization}. Thus, the geometric generic fiber of $\cX$ admits no integral decomposition of the diagonal, and in particular it is not stably rational over $\overline{\mathrm{Frac}(R)}$.
\end{prop}
   \begin{proof}
       By construction, $\cX$ is projective and flat over $\Spec(R)$. We note that $h'$ in~\eqref{eqn:Zcirc} has degree 4 in variables $T_{g(\mu)}$ and $T_l$ since $J$ in~\eqref{eqn:Ftau} is a general quartic polynomial. Let $\cX_0$ be the special fiber. In $S^\circ\times\bP^r$, it is given by 
       $$
      \left\{V^2\left( \sum_{\mu\in\mathcal{P}}\mu T
   _{g(\mu)}^2
   + \sum_{j=1}^{n}\mu_jT_j^2
+ (f_1-t^2f_2)V^2\right)=0\right\}\subset S^\circ\times\bP^r,
       $$
where we put coordinates $[V:T_{g(\mu)}:T_l]$ on $\bP^r$.
In particular $\cX_0$ contains an irreducible component $Y$ with multiplicity one, and the projection $Y\to S$ is proper and surjective with generic fiber  given by 
$$
Q=\left\{\sum_{\mu\in\mathcal{P}}\mu T
   _{g(\mu)}^2
   + \sum_{j=1}^{n}\mu_jT_j^2
+ (f_1-t^2f_2)V^2=0\right\}\subset\bP^r_{L(S)}.
$$
Note that $Q$ is exactly the quadric as in Proposition~\ref{prop:quadric}. The assumptions (1) and (2) in Theorem~\ref{prop:specialization} are satisfied by Proposition~\ref{prop:quadric}. From the equations, one checks that $Q$ is smooth on the affine chart $\{V\ne 0\}$. Thus 
$$
(Y\cap\Sing(\cX_0))_{L(S)}\subset Q\cap\{V=0\} 
$$
and the assumption (3) of Theorem~\ref{prop:specialization} is also satisfied. 
   \end{proof}

\begin{theo}\label{thm:mainreapt}
Let $k$ be an uncountable algebraically closed field of characteristic different from 2, and let $X\subset\bP^N$ be a very general quartic hypersurface with 
$$
2^a3^b-1\leq N\leq 2^{a+2b}+a+2b-2,\quad \text{for some } a\geq0,b\geq 2.
$$   
Then $X$ does not admit an integral decomposition of the diagonal, and in particular it is not stably rational over $k$.
\end{theo}
\begin{proof}
By \cite[Proposition 2.1]{torsionorderCL} it is enough to find one smooth quartic hypersurface  $X_0\subset \mathbb P^N$ defined over an algebraically closed field $\Omega\supset k$, such that $X_0$ has no integral decomposition of the diagonal.

Let $R=\overline{k(t)}[[\tau]]$, let $\Omega=\overline{\mathrm{Frac}(R)}$, and  let
$$
\cZ=\{F_\tau=0\}\subset\bP^N,
$$ 
be 
the smooth quartic constructed in \eqref{eqn:cZ}.  Let $\cX\to\Spec(R)$ be as in Proposition \ref{prop:degen}: its geometric generic fiber 
has no integral decomposition of the diagonal ({\it loc. cit.}), and it is birational to
 $\cZ_{\Omega}$ by construction.   By Lemma \ref{lemm:singdod}, $\cZ_{\Omega}$ has no integral decomposition of the diagonal, so that taking $X_0=\cZ_{\Omega}$ finishes the proof.

\end{proof}

\begin{coro}\label{coro:quartic}
    Let $k$ be an uncountable algebraically closed field with $\mathrm{char}(k)\ne 2$. Then a very general quartic hypersurface over $k$ of any positive dimension does not admit a decomposition of the diagonal, and in particular it is not stably rational.  
\end{coro}
\begin{proof}
    The cases of quartics of dimensions $\geq 7$ follow from Theorem~\ref{thm:mainreapt}. Taking
    \[
    (a,b) = (0,2),\: (1,2),\: (0,3),\: (1,3),
    \]
    we obtain the result in the intervals of (quartic) dimensions
    \[
    [7,17],\: [16,34],\: [25, 67],\: [52, 132]. 
    \]
   Taking $a=0$ and $b \ge 4$, we obtain overlapping intervals 
    \[
    [3^b -2, 4^b +2b-3] .
    \]
    The cases of dimensions $\le 6$ follow from the previous work of \cite{ColliotThelenePirutka, SchreiederSmallSlopes, Nicaseottem,quart5sch,OttemQuarticSixfolds}, as 
discussed in the introduction. 
\end{proof}
\begin{rema}\label{rema:higher}
    By setting 
$$
 F_\tau
= W_1^d\left(\sum_{\mu\in\mathcal{P}}
   \left(\tau^{2(r_0-\omega_\mu)}F_\mu\right)
   + \tau^{2(r_0-1)}\sum_{j=1}^{n}\left(s_j^2Q_j\right)
+ \tau^{2r_0}F\right)
   + \tau^{2r_0+1}J
 $$   
 for $d\in\bZ_{\geq 0}$, and $J$ a very general homogeneous polynomial of degree $d+4$, the same argument above shows that a very general hypersurface of degree at least 4 and dimension at least 7 is not stably rational. This recovers the result of \cite{HotchkissDegreeFive}. 
\end{rema}

    \bibliography{quartic}
	\bibliographystyle{alpha}
	
\end{document}